\documentclass[11pt]{article}
\usepackage{latexsym, amssymb, amsmath, amsthm}

\usepackage[colorlinks=true,citecolor=blue,linkcolor=blue,urlcolor=blue,bookmarks,bookmarksopen,bookmarksdepth=2,backref=page]{hyperref}

\def \F {{\mathbb F}}
\def \Z {{\mathbb Z}}
\def \V {{\mathbb V}}

\def \T {{\rm Tr}}

\def \+ {\oplus}
\newtheorem{theorem}{Theorem}

\newtheorem{lemma}{Lemma}
\newtheorem{proposition}{Proposition}
\newtheorem{remark}{Remark}

\newtheorem{corollary}{Corollary}

\begin{document}

\begin{center}

{\Large \bf Bent and $\Z_{2^k}$-bent functions from spread-like partitions}\\[1em]
Wilfried Meidl, Isabel Pirsic \\[.7em]
\end{center}
Johann Radon Institute for Computational and Applied Mathematics, \\
Austrian Academy of Sciences, Altenbergerstrasse 69, 4040-Linz, Austria.\\
 e-mail: meidlwilfried@gmail.com; isa.pirsic@gmail.com \\

\begin{abstract}
Bent functions from a vector space $\V_n$ over $\F_2$ of even dimension $n=2m$ into the cyclic group $\Z_{2^k}$, or equivalently, relative
difference sets in $\V_n\times\Z_{2^k}$ with forbidden subgroup $\Z_{2^k}$, can be obtained from spreads of $\V_n$ for
any $k\le n/2$. In this article, existence and construction of bent functions from $\V_n$ to $\Z_{2^k}$, which do not come from the spread
construction is investigated. A construction of bent functions from $\V_n$ into $\Z_{2^k}$, $k\le n/6$, (and more generally, into any abelian
group of order $2^k$) is obtained from partitions of $\F_{2^m}\times\F_{2^m}$, which can be seen as a generalization of the Desarguesian
spread. As for the spreads, the union of a certain fixed number of sets of these partitions is always the support of a Boolean bent function.
\end{abstract}

\section{Introduction}

Let $(A,+_A)$, $(B,+_B)$ be finite abelian groups. A function $f$ from $A$ to $B$ is called a {\it bent function} if
\begin{equation}
\label{CS}
|\sum_{x\in A}\chi(x,f(x))| = \sqrt{|A|}
\end{equation}
for every character $\chi$ of $A\times B$ which is nontrivial on $B$. Alternatively, $f:A\rightarrow B$ is bent if
and only if for all nonzero $a\in A$ the function $D_af(x) = f(x +_A a) -_B f(x)$ is balanced, i.e., every value in $B$ is taken on the same
number $|A|/|B|$ times. The graph of $f$, $G=\{(x,f(x))\,:\,x\in A\}$, is then a relative difference set in $A\times B$ relative to $B$,
see \cite{p04}. For background on relative difference sets we refer to \cite{p96}.

In the classical case, $A = \V_n$ and $B = \V_m$ are elementary abelian $2$-groups, i.e., they are vector spaces of dimension $n$ and
$m$ respectively over the prime field $\F_2$.
In this case the character sum in $(\ref{CS})$, called Walsh transform of $f$ at $(a,b)\in\V_m^*\times\V_n$, is of the form
\[ \mathcal{W}_f(a,b) = \sum_{x\in \V_n}(-1)^{\langle a,f(x)\rangle_m \oplus{+} \langle b,x\rangle_n}, \]
where $\langle , \rangle_k$ denotes an inner product in $\V_k$. If $\V_k = \F_2^k$ we may use the conventional dot product, the standard inner
product in $\F_{2^k}$, the finite field of order $2^k$, is $\langle u,v\rangle_k = \T_k(uv)$, the absolute trace of $uv$.
A function $f:\V_n\rightarrow\V_m$ is bent, if $m > 1$ also called {\it vectorial bent}, if and only if $|\mathcal{W}_f(a,b)| = 2^{n/2}$ for
all nonzero $a\in\V_m$ and $b\in\V_n$. As is well known, $n$ must then be even and $m$ can be at most $n/2$. Throughout the article,
$n=2m$ shall always be an even integer.

For Boolean bent functions, i.e., bent functions $f$ from $\V_n$ to $\F_2$, the dual $f^*$ is the Boolean function defined by
$\mathcal{W}_f(1,b) = \mathcal{W}_f(b) = 2^{n/2}(-1)^{f^*(b)}$, which is always a Boolean bent function as well.
There are many examples and constructions of Boolean bent functions in the literature. Even several classes of bent functions from $\V_n$
to $\V_{n/2}$ are known, such as Maiorana-McFarland functions, spread bent functions, Dillons $H$-class, see \cite{dillon}, and Kasami
bent functions, cf.\cite{c11}. For background on Boolean and vectorial (Boolean) bent functions we also refer to \cite{mesbook}.

In this article we are particularly interested in bent functions $f$ from $\V_n$ to the cyclic group $\Z_{2^k}$, hence by $(\ref{CS})$, in functions $f$ for which
\begin{equation}
\label{H_f}
\mathcal{H}_f(a,b) = \sum_{x\in\V_n}\zeta_{2^k}^{af(x)}(-1)^{\langle b,x\rangle},
\end{equation}
where $\zeta_{2^k}$ is a complex primitive $2^k$th root of unity, has absolute value $2^{n/2}$ for all nonzero $a\in\Z_{2^k}$ and $b\in\V_n$. Again such
functions can only exist for $m\le n/2$, \cite{n,sb}. %Hence in this article $n=2m$ will always be even.

We remark that functions $f:\V_n\rightarrow\Z_{2^k}$ satisfying the weaker condition that $|\mathcal{H}_f(1,b)| = 2^{n/2}$ for all $b\in\V_n$ are referred to
as generalized bent functions. They have been intensively studied in many papers, see \cite{hmp,mms,mms1,m18,mp,mtqwwf,txqf}.
If not also bent, generalized bent functions do not correspond to relative difference sets. However, as easily observed, $f:\V_n\rightarrow\Z_{2^k}$ is bent if
and only if $2^tf$ is a generalized bent function for every $t$, $0\le t\le k-1$.

As is well known, see also Section \ref{prel}, with the spread construction one can obtain bent functions from $\V_n$ to any abelian group of order $2^{n/2}$,
in particular also bent functions from $\V_n$ to $\Z_{2^{n/2}}$ (and as their projections, bent functions into cyclic groups of smaller order).
Unlike in the case of bent functions between elementary abelian groups, it seems difficult to find other examples for bent functions from $\V_n$ to $\Z_{2^k}$
when $k\ge 3$. In \cite{m18} it is observed that one example of a bent function from $\V_n$ to $\Z_{8}$, which does not come from (partial) spreads, can be
obtained from a secondary construction of Boolean bent functions in \cite{mes}.
As also pointed out in \cite{mp}, finding bent functions into the cyclic group which are not related to spreads is an interesting problem.

In Section \ref{prel} we revisit the construction of bent functions via spreads, in particular we point to constructing vectorial bent functions from $\V_n$ to $\V_m$,
and more general of bent functions from $\V_n$ into arbitrary abelian groups $G$, $|G| = 2^k$, obtained from partial, but not complete spreads.
In Section \ref{nospread}, a construction of bent functions from $\V_n = \F_{2^m}\times\F_{2^m}$ into $\Z_{2^k}$ is given in polynomial form. With an argument
using algebraic degree we show that this construction yields bent functions into $\Z_{2^k}$, $k\le m/3$, which cannot be obtained from (partial) spreads.
In Section \ref{partition}, analysing partitions of $\V_n$ coming with the functions of Section \ref{nospread}, we obtain large classes of bent functions
(Boolean, vectorial, and into the cyclic group), which have similar properties as spread functions, in fact can be seen as a generalization of spread functions obtained
from the Desarguesian spread. In particular, the union of a certain fixed number of sets of these partitions is always the support of a Boolean bent function.

%\section{Relative difference sets from partial spreads}
\section{Preliminaries}
\label{prel}

Recall that a partial spread $\mathcal{S}$ of $\V_n$, $n=2m$, is a set of $m$-dimensional subspaces of $\V_n$ which pairwise intersect trivially.
If $|\mathcal{S}| = 2^m+1$, hence every nonzero element of $\V_n $ is in exactly one of those subspaces, then $\mathcal{S}$ is called a (complete) spread.
The standard example is the Desarguesian spread, which has for $\V_n = \F_{2^m}\times\F_{2^m}$ the representation $\mathcal{S} = \{U, U_s \,:\,s\in\F_{2^m}\}$,
with $U = \{(0,y)\,:\,y\in\F_{2^m}\}$ and for $s\in\F_{2^m}$, $U_s = \{(x,sx)\,:\,x\in\F_{2^m}\}$.

As Dillon showed in his thesis \cite{dillon}, a Boolean bent function $f$ is obtained by choosing the nonzero elements of $2^{m-1}$ subspaces of $\mathcal{S}$ as the
support of $f$ (PS$^-$ bent functions), or the elements of $2^{m-1}+1$ subspaces of $\mathcal{S}$ for the support of $f$ (PS$^+$ bent functions).
Clearly, every Boolean bent function which is constant on the subspaces of a (partial) spread is of this form.
As is also well known, with a complete spread $\mathcal{S}$ one gets a vectorial bent function from $\V_n$ to $\V_m$ by mapping to every nonzero $z\in \V_m$
the nonzero elements of exactly one subspace of $\mathcal{S}$,  the remaining two subspaces are mapped to $0$.
From projections of a bent function from  $\V_n$ to $\V_m$ one can get vectorial bent functions from $\V_n$ to $\V_k$ for any $1\le k\le n/2$. We remark that
the spread construction can also be applied to elementary abelian $p$ groups, $p$ odd, see for instance \cite{ll} for PS$^-$ and PS$^+$ bent functions from
$\F_p^{n}$ to $\F_p$. The above extreme cases, Boolean bent functions for which solely $2^{m-1}$, respectively $2^{m-1}+1$ subspaces are required, and the
construction of vectorial bent functions from $\V_n$ to $\V_m$ with complete spreads are mostly considered in the literature. However, the spread construction can also be applied with
not complete spreads (with more elements) to obtain vectorial bent functions. The resulting functions are then in general not projections of a vectorial spread bent
function into $\V_m$. The procedure is standard, but as far as we know, the argument has not been given explicitly in the literature on bent functions, hence we give it below.
%we find it worth to explicitly point to this fact at this position.
We consider the most general situation, i.e., functions from $\V_n$ to $B$, where $B$ is any abelian group of order $2^k$, $1\le k\le m$.
The bentness condition is then $|\mathcal{T}_f(b)| = 2^m$ for all $b\in\V_n$ and all nontrivial
characters $\chi$ of $B$, where, for short, denoting the inner product $\langle u,v\rangle$ in $\V_n$ by $\langle u,v\rangle = u\cdot v$,
\[ \mathcal{T}_f(b) = \sum_{x\in \V_n}\chi(f(x))(-1)^{b\cdot x}. \]
Let $\mathcal{S} = \{U_j,\,1\le j\le (2^k-1)2^{m-k}\}$ be a partial spread of $\V_n$ and $B$ an abelian group of order $2^k$.
Define $f:\V_n\rightarrow B$ as follows: \\[.3em]
Construction I.
\begin{itemize}
\item[-] Every nonzero element $\gamma$ of $B$ has as preimage the union of exactly $2^{m-k}$ elements of $\mathcal{S}$
except from $0\in\V_n$, i.e., $f^{-1}(\gamma) = \bigcup_{i=1}^{2^{m-k}}U^*_{\gamma,i}$, where $U^* = U\setminus\{0\}$.
\item[-] All other elements are mapped to $0\in B$, i.e., $f^{-1}(0) = \V_n\setminus\bigcup_jU^*_j$.
\end{itemize}
First observe that
\begin{align*}
\mathcal{T}_f(0) & = \sum_{\gamma\in B^*}\sum_{x\in U^*_{\gamma,i}\atop 1\le i\le 2^{m-k}}\chi(\gamma) + 2^{2m} - \sum_{x\in\bigcup_j U_j^*}1 \\
& = (-1)2^{m-k}(2^m-1) + 2^{2m} - (2^k-1)2^{m-k}(2^m-1) = 2^m.
\end{align*}
If $b\ne 0$, then
\begin{equation}
\label{T=}
\mathcal{T}_f(b) %= \sum_{\gamma\in G^*}\sum_{x\in U^*_{\gamma,i}\atop 1\le i\le 2^{m-k}}\chi(\gamma)(-1)^{b\cdot x} + \sum_{x\in D_0}(-1)^{b\cdot x}
= \sum_{\gamma\in B^*}\sum_{x\in U^*_{\gamma,i}}\chi(\gamma)(-1)^{b\cdot x} - \sum_{x\in\bigcup_j U_j^*}(-1)^{b\cdot x}.
\end{equation}
We use that for every nonzero $b\in\V_n$ we have at most one $U_j\in\mathcal{S}$ for which $b\cdot x = 0$ for all $x\in U_j$. \\[.3em]
{\it Case 1:} There is no $U_j\in\mathcal{S}$ such that $b\cdot x = 0$ for all $x\in U_j$. \\
Then we have
\[ \mathcal{T}_f(b) = \sum_{\gamma\in B^*}\chi(\gamma)2^{m-k}(-1) - (2^k-1)2^{m-k}(-1) = 2^m. \]
{\it Case 2:} $b\cdot x = 0$ for all $x\in U_{\gamma_0,i}$ for (exactly) one $\gamma_0\in B^*$ and one $i\in\{1,\ldots,2^{m-k}\}$. \\
In this case,
\begin{align*}
\mathcal{T}_f(b) & = \sum_{\gamma\in B^*}\chi(\gamma)2^{m-k}(-1) + 2^m\chi(\gamma_0) - ((2^k-1)(-1)2^{m-k} + 1 + 2^m-1) \\
& = 2^{m-k} + 2^m\chi(\gamma_0) + (2^k-1)2^{m-k} - 2^m = 2^m\chi(\gamma_0).
\end{align*}

Similarly, one shows that for a partial spread with (at least) $2^m-2^{m-k}+1$ subspaces, the function $g:\V_n\rightarrow B$ defined
as follows is a bent function: \\[.3em]
Construction II.
\begin{itemize}
\item[-]  For an element $\tilde{\gamma} \in B^*$ we have $g^{-1}(\tilde{\gamma}) = \bigcup_{i=1}^{2^{m-k}+1}U_{\tilde{\gamma},i}$,
i.e., $\tilde{\gamma}$ has the union of $2^{m-k}+1$ elements of a partial spread as preimage (note that also $f(0) = \tilde{\gamma}$),
\item[-] if $\gamma \in B^*$, $\gamma \ne \tilde{\gamma}$, then $g^{-1}(\gamma) = \bigcup_{i=1}^{2^{m-k}}U^*_{\gamma,i}$,
i.e., the preimage of $\gamma$ consists of the nonzero elements of $2^{m-k}$ elements of a partial spread,
\item[-] the remaining elements are mapped to $0$.
\end{itemize}

If $k=1$, hence $B = \F_2$, then the functions in Construction I and II are conventional PS$^-$ and PS$^+$ bent functions, respectively.
Hence one may see the vectorial partial spread functions from $\V_n$ to $\V_k$ obtained by Construction I and II with $B=\V_k$ as {\it vectorial PS$^-$
and PS$^+$ bent functions}. Note that for Construction I respectively II one needs partial spreads $\mathcal{S}$ with at least $2^m-2^{m-k}$
respectively $2^m-2^{m-k}+1$ elements.
%NOW: Check literature about sizes of not extendable partial spreads. I would at least say, for $k=m$,
%if we have $2^m-1$ subspaces, then we must be in a complete spread.

In this article we are interested in bent functions from $\V_n$ to the cyclic group $\Z_{2^k}$.
Canonical examples one obtains with the partial spread construction with $B = \Z_{2^k}$ for all $k\le n/2$,
which we will also call partial spread bent functions from $\V_n$ to $\Z_{2^k}$.
In the following sections we investigate existence and construction of bent functions from $\V_n$ to $\Z_{2^k}$ which do not come from partial spreads.

\section{$\Z_{2^k}$-bent functions not obtained from spreads}
\label{nospread}

As we have to distinguish addition in different structures, we denote the addition in the complex numbers and in the ring $\Z_{2^k}$ by $+$,
the addition in the elementary abelian groups $\F_2$, $\V_n$ and $\F_{2^m}$ is denoted by $\oplus$.

Let $f$ be a function from $\V_n$ to $\Z_{2^k}$, then we can write $f$ as
\begin{equation}
\label{fform}
f(x) = a_0(x) + 2a_1(x) + \cdots + 2^{k-1}a_{k-1}(x)
\end{equation}
for uniquely determined Boolean functions $a_j$, $0\le j\le k-1$, from $\V_n$ to $\F_2$.

Recall that a function $f:\V_n\rightarrow\Z_{2^k}$ is called generalized bent if for $\mathcal{H}_f$ given as in $(\ref{H_f})$ we have
$|\mathcal{H}_f(1,u)| = 2^{n/2}$ for all $u\in\V_n$. Naturally, bent functions from $\V_n$ to $\Z_{2^k}$ are generalized bent. The
converse does not hold, but we have the following obvious lemma, see \cite{hmp}.
\begin{lemma}
\label{clear}
A function $f:\V_n\rightarrow\Z_{2^k}$ is bent if and only if $2^tf$ is generalized bent for all $t$, $0\le t\le k-1$.
\end{lemma}
Generalized bent functions are intensively studied in the literature, see \cite{hmp,mms,mms1,m18,mp,mtqwwf,txqf}.
A comprehensive characterization of generalized bent functions via partitions has been given in \cite{mtqwwf}.
In \cite{hmp}, another characterization of generalized bent functions has been given via properties of the Boolean functions
$a_i$ in $(\ref{fform})$.
\begin{proposition}
\label{enes}
Let $f:\V_n\rightarrow\Z_{2^k}$ be given as $f(x) = a_0(x) + 2a_1(x) + \cdots + 2^{k-1}a_{k-1}(x)$ for some Boolean functions
$a_j$, $0\le j\le k-1$, and let $\mathcal{A}$ be the affine space of Boolean functions,
\[ \mathcal{A} = a_{k-1} \oplus \langle a_{k-2},\ldots, a_0\rangle. \]
Then $f$ is generalized bent if and only if all functions in $\mathcal{A}$ are Boolean bent functions, and for any three functions
$b_0,b_1,b_2 \in \mathcal{A}$ we have
\begin{equation}
\label{bbb*}
(b_0\oplus b_1\oplus b_2)^* = b_0^*\oplus b_1^*\oplus b_2^*.
\end{equation}
\end{proposition}

One of the main objectives in this article is to construct and to understand infinite classes of bent functions from $\V_n$ to the cyclic group
$\Z_{2^k}$, or equivalently of relative difference sets in $\V_n\times\Z_{2^k}$ relative to $\Z_{2^k}$, which do not come from partial spreads.
We therefore need a tool to distinguish such bent functions from partial spread bent functions.
As in our analysis of bent functions from $\V_n$ to $\Z_{2^k}$, Boolean bent functions play a major role, for better understanding, in the remainder
of the article we will denote bent functions from $\V_n$ to $\Z_{2^k}$ as $\Z_{2^k}$-bent functions. Bent functions from $\V_n$ to $\F_2$ will
be referred to as Boolean bent functions or simply as bent functions.

Recall that two Boolean functions $f,g$ from $\V_n$ to $\F_2$ are called {\it extended affine equivalent (EA-equivalent)}
if $g(x) = f(\mathcal{L}(x)+b) + a(x)$ for some linear coordinate transformation $\mathcal{L}$ on $\V_n$, an element
$b\in\V_n$ and an affine map $a(x)$ from $\V_n$ to $\F_2$. If $a(x)$ is the zero-map, then $f$ and $g$ are called
{\it afflne equivalent}, if additionally $b = 0$, then $f$ and $g$ are called {\it linear equivalent}.
\begin{lemma}
\label{isPS}
Let $f:\V_n\rightarrow\Z_{2^k}$ with $f(x) = a_0(x) + 2a_1(x) + \cdots + 2^{k-1}a_{k-1}(x)$ be a partial spread $\Z_{2^k}$-bent function.
Then all Boolean functions $a_i$, $0\le i\le k-1$, are Boolean partial spread bent functions, all of algebraic degree $m=n/2$.
\end{lemma}
{\it Proof.}
Since $f$ is $\Z_{2^k}$-bent, by Lemma \ref{clear}, for every $t$, $0\le t\le k-1$, every Boolean function in $\mathcal{A}_t = a_{k-t-1} \oplus \langle a_{k-t-2},\ldots, a_0\rangle$
is bent. In particular, $a_i$, $0\le i\le k-1$, is a Boolean bent function.
Let $U \in \mathcal{S}$, and $x_1,x_2\in U^*$. If $a_j(x_1) \ne a_j(x_2)$, then $f:\V_n\rightarrow\Z_{2^k}$ cannot be constant on $U^*$.
Hence $a_j$ is a Boolean bent function which is constant on the nonzero elements of $U$ for every $U\in\mathcal{S}$, consequently a
Boolean partial spread bent function. By \cite[p.96]{dillon}, every Boolean partial spread bent function on a spread with more than $2^{m-1}$
subspaces has algebraic degree $m$.
\hfill$\Box$
\begin{remark}
With the same argument, all bent functions in $\mathcal{A}_t$, $0\le t\le k-1$, are partial spread bent functions
of algebraic degree $m$.
\end{remark}

In \cite{mes}, Mesnager presented several examples of Boolean bent functions $b_0,b_1,b_2$ satisfying
$(b_0\oplus b_1\oplus b_2)^* = b_0^*\oplus b_1^*\oplus b_2^*$. The objective in \cite{mes} is to use those functions in a secondary construction of Boolean
bent functions due to Carlet \cite{c06}. As shown in \cite{m18} this secondary construction is equivalent to constructing
generalized bent functions to $\Z_{2^3}$. As observed in \cite{m18,mp}, one of the examples, using Maiorana-McFarland functions, potentially
yields $\Z_{2^3}$-bent functions from $\V_n = \F_{2^m}\times\F_{2^m}$ to $\Z_{2^3}$.
Recall that for a permutation $\pi$ of $\F_{2^m}$, the Boolean function $h(x,y) = \T_m(x\pi(y))$ from $\F_{2^m}\times\F_{2^m}$ to $\F_2$ is a bent
function belonging to the Maiorana-McFarland class. Noting that for the Maiorana-McFarland bent function $b(x,y) = \T_m(\beta xy^d)$, $\gcd(2^m-1,d) = 1$,
from $\F_{2^m}\times\F_{2^m}$ to $\F_2$, we have $b^*(x,y) = \T_m(\beta^{-e}x^ey)$, $ed\equiv 1 \bmod 2^m-1$, we arrive at the following observation
which we here state as a lemma. By convention, all powers of $0$ are equal to $0$ (including powers with 
negative exponents). 
\begin{lemma}\cite{mes}
\label{mes}
Let $d,e$ be integers such that $\gcd(2^m-1,d) = 1$ and $ed \equiv 1 \bmod 2^m-1$, and suppose that
$\beta_0,\beta_1,\beta_2 \in\F_{2^m}$ satisfy
\begin{equation}
\label{bebebe}
(\beta_0 \oplus \beta_1 \oplus \beta_2)^{-e} = \beta_0^{-e} \oplus \beta_1^{-e} \oplus \beta_2^{-e}.
\end{equation}
Then the Boolean bent functions $b_i(x) = \T_m(\beta_i xy^d)$, $i = 0,1,2$, satisfy
$(b_0\oplus b_1\oplus b_2)^* = b_0^*\oplus b_1^*\oplus b_2^*$.
\end{lemma}
As easily observed, $f(x) = b_0(x) + 2(b_0 \oplus b_1)(x) + 4(b_0\oplus b_2)(x)$ is then generalized bent, and since $b_i \oplus b_j$, $0\le i < j \le 2$,
is bent, $f$ is even $\Z_8$-bent. For the details we refer to \cite{m18}.

Trivially, with $-e = -1$, and more general $-e \equiv 2^v\bmod 2^m-1$, equation $(\ref{bebebe})$ is satisfed
for all choices of $\beta_0,\beta_1,\beta_2 \in\F_{2^m}$.
In \cite[Table 1]{mes}, for $4\le n\le 8$ some exponents $e$ for which there exist $\beta_0,\beta_1,\beta_2$ such that $(\ref{bebebe})$ is satisfied are
listed.

In fact we have recalculated the values and slightly expanded the table, also
completing the table with
some additional entries not present in the original. We find that the
examples are so prolific that it is in all but two cases more efficient to list the
complements, i.e., the coset leaders of the cyclotomic classes that do
\underline{not} fulfill
the condition (also we omit the class led by $1$).

\begin{table}[ht]\begin{center}\begin{tabular}{|c|c|l|}
\hline
 $n$ & $2^n-1$ & Cyclotomic class leaders \\
\hline\hline
 4 & 15 & 3 \\ \hline
 5 & 31 & only 15 \underline{fulfills}  the cond.\\ \hline
 6 & 63 & 15 \hspace{.52cm} (\underline{not} fulfilling \eqref{bebebe}) \\ \hline
 7 & 127 & 3, 5, 9, 15, 27, 43, 63 \underline{fulfill} the cond. \\ \hline
 8 & 255 & 27, 63, 111 \\ \hline
 9 & 511 & 15, 29, 39, 51, 53, 79, 85, 95, 123, 127, 191, 239\\\hline
 10 & 1023 & 111, 171, 255, 447 \\\hline
\end{tabular}\caption{List of classes of exponents $e$
\underline{not} fulfilling \eqref{bebebe}}
\end{center}\end{table}

With the next lemma we obtain large sets of elements of which every three satisfy the condition $(\ref{bebebe})$ for some fixed (nontrivial) $e$. In fact even stronger
conditions are satisfied, which we will require to construct $\Z_{2^k}$-bent functions for $k > 3$.
\begin{lemma}
\label{beholds}
Let $m, j$ be integers such that and $\gcd(2^m-1,2^j-1) = 2^k-1$, and let $e = 2^m-2^j-2$.
Then for any sum $\sum_i \beta_i$ of elements of $\F_{2^k}$ we have $\sum\beta_i^{-e} = (\sum\beta_i)^{-e}$.
In particular, for any $\beta_0,\beta_1,\beta_2 \in\F_{2^k}$ we have
$(\beta_0 \oplus \beta_1 \oplus \beta_2)^{-e} = \beta_0^{-e} \oplus \beta_1^{-e} \oplus \beta_2^{-e}$.
\end{lemma}
{\it Proof.}
As $k$ divides $j$, for $\beta\in \F_{2^k}$ we have $\beta^{-e} = \beta^{2^j+1} = \beta^2$, hence $\sum\beta_i^{-e} = (\sum\beta_i)^{-e}$
if $\beta_i\in\F_{2^k}$. In particular we have $(\beta_0 \oplus \beta_1 \oplus \beta_2)^{-e} = \beta_0^{-e} \oplus \beta_1^{-e} \oplus \beta_2^{-e}$
for $\beta_0,\beta_1,\beta_2 \in\F_{2^k}$.
\hfill$\Box$\\[.5em]
Before we show the main results of this section, we apply Lemma \ref{beholds} to the secondary construction of bent functions in \cite{c06}, and obtain
with Theorem 14 in \cite{mes} the following
\begin{corollary}
\label{Carletsec}
Let $m, j$ be integers such that $\gcd(2^m-1,2^j+1) = 1$ and $\gcd(2^m-1,2^j-1) = 2^k-1$, let $e = 2^m-2^j-2$, and let $d$ be the inverse of $e$ modulo $2^m-1$.
Then for any $\beta_0, \beta_1, \beta_2 \in\F_{2^k}^*$ the Boolean function
\begin{align*}
g(x) & = \T_m(\beta_0^{-e}x^ey)\T_m(\beta_1^{-e}x^ey) \oplus \T_m(\beta_0^{-e}x^ey)\T_m(\beta_2^{-e}x^ey) \\
& \oplus \T_m(\beta_1^{-e}x^ey)\T_m(\beta_2^{-e}x^ey)
\end{align*}
is bent. Its dual is
\begin{align*}
g^*(x) & = \T_m(\beta_0xy^d)\T_m(\beta_1xy^d) \oplus \T_m(\beta_0xy^d)\T_m(\beta_2xy^d) \\
& \oplus \T_m(\beta_1xy^d)\T_m(\beta_2xy^d).
\end{align*}
\end{corollary}
We now come back to our primary objective, the construction of $\Z_{2^k}$-bent functions, and state the main results of this section.
\begin{theorem}
\label{main}
Let $m, j$ be integers such that $\gcd(2^m-1,2^j+1) = 1$ and $\gcd(2^m-1,2^j-1) = 2^k-1$, let $e = 2^m-2^j-2$, and let
$d$ be the inverse of $e$ modulo $2^m-1$. Then for a basis $\{\alpha_0,\alpha_1,\ldots,\alpha_{k-1}\}$ of $\F_{2^k}$ over $\F_2$,
the functions $f_1$ and $f_2$ given as
\begin{equation}
\label{f1f2}
f_1(x) = \sum_{i=0}^{k-1}\T_m(\alpha_i xy^d)2^i, \qquad f_2(x) = \sum_{i=0}^{k-1}\T_m(\alpha_i^{-e} x^ey)2^i
\end{equation}
are $\Z_{2^k}$-bent functions from $\F_{2^m}\times\F_{2^m}$ to $\Z_{2^k}$.
\end{theorem}
{\it Proof.}
By Lemma \ref{clear}, a function $f:\V_n\rightarrow\Z_{2^k}$ given as $f(x) = a_{k-1}(x)2^{k-1} + a_{k-2}(x)2^{k-2} + \cdots + a_0(x)$ is
$\Z_{2^k}$-bent if and only if $2^tf(x) = a_{k-t-1}(x)2^{k-1} + \cdots + a_0(x)2^t$ is generalized bent for all $0\le t\le k-1$.
With Proposition \ref{enes}, the function $2^tf(x)$ is generalized bent if and only if every Boolean function in the affine space
$\mathcal{A}_t = a_{k-t-1} \oplus \langle a_{k-t-2},\ldots,a_0\rangle$ is bent, and for each three Boolean bent functions $b_0, b_1, b_2 \in\mathcal{A}_t$
we have $(b_0\oplus b_1\oplus b_2)^* = b_0^*\oplus b_1^*\oplus b_2^*$. (Note that then also $b_0\oplus b_1\oplus b_2$ is a bent function in $\mathcal{A}_t$.)

Observe that for the function $f_1$, every Boolean function in $\mathcal{A}_t$, $0\le t\le k-1$, is of the form $b(x) = \T_m(\beta xy^d)$, where
$\beta$ is an element in $\F_{2^k}^*$, thus $b(x)$ is bent. Let now $b_0, b_1, b_2 \in\mathcal{A}_t$ be given as $b_i(x) = \T_m(\beta_i xy^d)$, $i=0,1,2$.
By Lemma \ref{beholds},
equation $(\ref{bebebe})$ is satisfied, and by Lemma \ref{mes} the function $f_1$ is $\Z_{2^k}$-bent.

For the function $f_2$, the Boolean functions in $\mathcal{A}_t$, $0\le t\le k-1$, are of the form $b(x) = \T_m(\tilde{\beta} x^ey)$, $\tilde{\beta} \in \F_{2^k}^*$.
Since $\gcd(-e,2^k-1) = 1$, hence $\tilde{\beta} = \beta^e$ for some $\beta \in\F_{2^k}$, for three functions $b_0,b_1,b_2 \in \mathcal{A}_t$ we can assume that
$b_i(x) = \T_m(\beta_i^{-e} x^ey)$, $\beta_i\in\F_{2^m}^*$, $i=0,1,2$. Then $b_0\oplus b_1\oplus b_2 = \T_m((\beta_0^{-e}\oplus \beta_1^{-e}\oplus \beta_2^{-e}) x^ey) =
\T_m((\beta_0\oplus \beta_1\oplus \beta_2)^{-e}x^ey)$. For the duals, $b_i^*(x) = \T_m(\beta_i xy^d)$, $i=0,1,2$, and
$(b_0\oplus b_1\oplus b_2)^*(x) = \T_m((\beta_0 \oplus \beta_1 \oplus \beta_2) xy^d)$, the condition $(\ref{bbb*})$ is satisfied. Consequently, by Proposition \ref{enes}, $f_2$ is
$\Z_{2^k}$-bent.
\hfill$\Box$
\begin{remark}
\label{rem2}
For $j=0$ the conditions of Theorem \ref{main} are satisfied. In this case, $k=m$ and one obtains a bent function from $\V_n$ to
$\Z_{2^m}$, $n = 2m$. As easily observed, $e = 2^m-3$ (hence $-e\equiv 2\bmod 2^m-1$ trivially satisfies 
condition $(\ref{bebebe})$), its inverse $d$ and $2^m-2$ belong to the same cyclotomic coset of $2$ modulo $2^m-1$,
i.e., $\T_m(\alpha xy^d) = \T_m(\alpha x^{2^v(2^m-2)}y)$ for some integer $v$. Hence $\T_m(\alpha xy^d)$ is
obtained form the classical Boolean partial spread bent function $\T_m(\alpha x^{2^m-2}y)$ by applying a 
linear coordinate transformation to the variable $x$. Therefore, $\T_m(\alpha xy^d)$ and $\T_m(\alpha x^{2^m-2}y)$
are affine equivalent. The functions $f_1,f_2$ in Theorem \ref{main} are then partial spread $\Z_{2^k}$-bent functions from $\F_{2^m}\times\F_{2^m}$
to $\Z_{2^m}$ as constructed in Section \ref{prel}. \\
Clearly, with $e=d=2^m-2$, the functions $f_1,f_2$ in $(\ref{f1f2})$ are partial spread $\Z_{2^k}$-bent functions, see also \cite[Corollary 17]{mes}.
\end{remark}
With Theorem \ref{main}, for $j>0$ on the other hand, we achieve our target and obtain $\Z_{2^k}$-bent functions from $\F_{2^m}\times\F_{2^m}$ to $\Z_{2^k}$
which provably do not come from partial spreads.
\begin{corollary}
Let $m$ and $j>0$ be integers such that $\gcd(2^m-1,2^j+1) = 1$ and $\gcd(2^m-1,2^j-1) = 2^k-1$, and let $e,d$, $\alpha_i$, $0\le i\le k-1$, be as in Theorem \ref{main}.
Then the functions $f_1,f_2$ in $(\ref{f1f2})$ are $\Z_{2^k}$-bent functions from $\F_{2^m}\times\F_{2^m}$ to $\Z_{2^k}$, which do not come from partial spreads.
\end{corollary}
{\it Proof.}
If $j>0$, then the binary weight of $e = 2^m-2^j-2$ is $m-2$. Hence the algebraic degree of $\T_m(\alpha_i^{-e} x^ey)$ in $f_2$ is $m-1$. By Lemma \ref{isPS},
$f_2$ is not a partial spread $\Z_{2^k}$-bent function. The bent functions $\T_m(\alpha_i xy^d)$ in $f_1$ are the duals of the bent functions $\T_m(\alpha_i^{-e} x^ey)$
in $f_2$, hence also not partial spread bent functions. Again with Lemma \ref{isPS}, $f_2$ is not a partial spread $\Z_{2^k}$-bent function.
\hfill$\Box$\\[.5em]
The largest possible $k$ in the construction of Theorem \ref{main}, besides from $k=m$, which as explained above solely yields the known spread functions,
is $k = m/3$. As easily seen, for every integer $m$ divisible by $3$ and $k=m/3$ we have $\gcd(2^k+1,2^m-1) = 1$. Hence the conditions of Theorem \ref{main}
are satisfied and we have the following corollary.
\begin{corollary}
Let $n=2m$ be divisible by $3$, then there exists a bent function from $\V_n$ to $\Z_{2^{m/3}}$ which is not obtained from a (partial) spread.
\end{corollary}
We expect that $\Z_{2^m}$-bent functions from $\V_{2m}$ to $\Z_{2^m}$ one only can obtain from (complete) spreads.

We finish this section with a remark on the vectorial bent function $F(x) = (\T_m(\alpha_0 xy^d),\T_m(\alpha_1 xy^d),\ldots,\T_m(\alpha_{k-1} xy^d))$
%\qquad f_2(x) = \sum_{i=0}^{k-1}\T_m(\alpha_i^{-e} x^ey)
from $\F_{2^m}\times\F_{2^m}$ to $\F_2^k$ associated with $f_1$, a projection onto $\F_{2^k}$ of the vectorial Maiorana-McFarland bent function
$\tilde{F}(x) = \T_m^n(xy^d)$ from $\F_{2^m}\times\F_{2^m}$ to $\F_{2^m}$. As shown in \cite[Theorem 2]{cmp18}, $\tilde{F}$ is a
{\it vectorial dual-bent function}, i.e., a vectorial bent function for which the set of the duals of all component functions again forms a vectorial bent
function (of the same dimension). The projection $F$ attached to our $\Z_{2^k}$-bent function satisfies the even stronger condition that
$(b\oplus \bar{b})^* = b^* \oplus \bar{b}^*$ for all component functions $b,\bar{b}$ of $F$.

\section{More bent and $\Z_{2^k}$-bent functions from the partitions}
\label{partition}

In this section we look at the partitions of $\F_{2^m}\times\F_{2^m}$ which we obtain from the sets of the preimages of the $\Z_{2^k}$-bent functions
$f_1$ and $f_2$ in Theorem \ref{main}. As we will see, these partitions share some properties with spreads, they can be seen as a generalization of the
Desarguesian spread. For simplicity we take $j=k$ a divisor of $m$.

For integers $m, e$ with $\gcd(2^m-1,e) = 1$ and an element $s\in\F_{2^m}$ define
\[ U_s :=\{(x,sx^{-e})\;:\;x\in \F_{2^m}\}, \; U_s^*=U_s\setminus\{0\}, \;\mbox{and}\; U = \{(0,y)\;:\;y\in\F_{2^m}\}. \]
Then $U$, $U_s^*$, $s\in\F_{2^m}$, form a partition of $\F_{2^m}\times\F_{2^m}$. Note that $U$, $U_s$, $s\in\F_{2^m}$,
are the subspaces of the Desarguesian spread if $-e\equiv 1\bmod 2^m-1$ (more general, if $-e\equiv 2^v\bmod 2^m-1$). Also
note that $U_s$ is not a subspace if we do not have $-e\equiv 2^v\bmod 2^m-1$ for some integer $v$.

Similarly, for integers $m, d$ with $\gcd(2^m-1,d) = 1$ and an element $s\in\F_{2^m}$ define
\[ V_s :=\{(x^{-d}s,x)\;:\;x\in \F_{2^m}\}, \; V_s^*=V_s\setminus\{0\}, \;\mbox{and}\; V = \{(x,0)\;:\;x\in\F_{2^m}\}. \]
Note that as above for the sets $U$ and $U_s$, if $-d\equiv 2^v\bmod 2^m-1$, then $V_s$ and $V$ are the subspaces of the Desarguesian spread.

For a divisor $k$ of $m$ and an element $\gamma$ of $\F_{2^k}$ let
\begin{equation}
\mathcal{A}(\gamma) = \bigcup_{s\in\F_{2^m}\atop\T^m_k(s)=\gamma}U_s^* \quad\mbox{and}\quad
\mathcal{B}(\gamma) = \bigcup_{s\in\F_{2^m}\atop\T^m_k(s)=\gamma}V_s^*.
\end{equation}
\begin{lemma}
\label{Gamma}
Let $k$ be a divisor of $m$, $\gcd(2^m-1,2^k+1) = 1$ and let $e=2^m-2^k-2$ and $d$ be the integer such that
$de \equiv 1\bmod 2^m-1$. For the $\Z_{2^k}$-bent functions $f_1$ and $f_2$ in Theorem \ref{main} let
\[ \Gamma_1 = \{A(i) = \{(x,y)\in\F_{2^m}\times \F_{2^m}\;:\;f_2(x,y) = i\}; i=0,\ldots2^k-1\} \]
be the partition of $\F_{2^m}\times \F_{2^m}$ obtained with the set of the preimages for $f_2$, and let
\[ \Gamma_2 = \{B(i) = \{(x,y)\in\F_{2^m}\times \F_{2^m}\;:\;f_1(x,y) = i\}; i=0,\ldots2^k-1\} \]
be the partition of $\F_{2^m}\times \F_{2^m}$ obtained with the set of the preimages for $f_1$. Then
\begin{align*}
\Gamma_1 & = \{\mathcal{A}(0)\cup U, \mathcal{A}(\gamma);\gamma\in\F_{2^m}^*\} \\
\Gamma_2 & = \{\mathcal{B}(0)\cup V, \mathcal{B}(\gamma);\gamma\in\F_{2^m}^*\},
\end{align*}
where $\mathcal{A}(0)\cup U = A(0)$ and $\mathcal{B}(0)\cup V = B(0)$.
\end{lemma}
{\it Proof.}
First observe that on $U^*_s = \{(x,sx^{2^k+1})\,:\,x\in\F_{2^m}^*\}$ the bent function $\T_m(\alpha^{-e}x^ey) = \T_m(\alpha^{-e}s)$ is constant, hence $f_2$
is constant on $U^*_s$. Let $s_1,s_2\in\F_{2^m}$, then $f_2$ takes on the same value on $U_{s_1}$ and $U_{s_2}$ if and only if $\T_m(\alpha_i^{-e} s_1) = \T_m(\alpha_i^{-e} s_2)$
for all $\alpha_i$ of the basis $\{\alpha_0,\ldots,\alpha_{k-1}\}$ of $\F_{2^k}$.
Hence we require
\begin{equation}
\label{s12}
\T_m(\alpha_i^{-e} s_1 \oplus \alpha_i^{-e} s_2) = \T_k(\alpha_i^{-e}\T_k^m(s_1\oplus s_2)) = 0\;\mbox{for all}\;i = 0, \ldots, k-1.
\end{equation}
Using that $\sum_j\alpha_{i_j}^{-e} = (\sum_j\alpha_{i_j})^{-e}$, we infer that $\{\alpha_i^{-e}; i=0,\ldots,k-1\}$ is a basis of $\F_{2^k}$ as well.
Hence $(\ref{s12})$ holds if and only if $\T_k^m(s_1\oplus s_2) = 0$. Consequently, $f_2$ is constant on
$\mathcal{A}(\gamma)$ and $\mathcal{A}(\gamma_1)$, $\mathcal{A}(\gamma_2)$ with $\gamma_1\ne\gamma_2$ are mapped to different constants.
Additionally, the elements of $U$, as well as the elements of $\mathcal{A}(0)$ are mapped to $0$. The same argument applies for $\Gamma_2$ with $f_1$.
\hfill$\Box$ \\[.5em]
Note that for $(x,y)\in U^*_s$ we have $\T_m(\alpha xy^d) = \T_m(\alpha x s^dx^{-ed}) = \T_m(\alpha s^d)$. Hence also $f_1$ is constant on $U^*_s$ and
certainly on $U$. In fact we can write $\Gamma_2$ also in terms of the sets $U$ and $\bigcup_{s\in\F_{2^m}\atop\T^m_k(s^d)=\gamma}U_s^*$.
Observe that the partitions $\Gamma_1$ and $\Gamma_2$, though looking similar, are different.
In the further we will use the representation of $\Gamma_2$ as in Lemma \ref{Gamma}.
\begin{remark}
In the special case $k=m$, the partition $\{U,\mathcal{A}(\gamma); \gamma\in\F_{2^m}\} = \{U,U^*_s; s \in\F_{2^m}\}$ with $U_s = \{(x,sx^2)\;:\;x\in\F_{2^m}\}$,
reduces to (a representation of) the Desarguesian spread partition of $\F_{2^m}\times\F_{2^m}$. The same applies to the partition
$\{V,\mathcal{B}(\gamma); \gamma\in\F_{2^m}\}$. As already observed in Remark \ref{rem2}, $f_1$ and $f_2$ are then spread $\Z_{2^m}$-bent functions.
As pointed out in Section \ref{prel}, many more can be obtained with the partition, by taking as preimage of every element in $\Z_{2^m}$ exactly one of these sets,
except for one (w.l.o.g. $0$), which has the union of two elements of the partition as preimage.
\end{remark}
To deduce a comparable result for the case $k<m$ we will need the following lemma.
\begin{lemma}
\label{aboutO}
Let $m,k$ be integers such that $k$ divides $m$ and $\gcd(2^m-1,2^k+1) = 1$, let $e = 2^m-2^k-2$ and $d$ be the inverse of $e$ modulo $2^m-1$, i.e.,
$de \equiv 1\bmod 2^m-1$. Let $u,v \in\F_{2^m}$ and $\gamma\in\F_{2^k}$.
\begin{itemize}
\item[(i)] If $v\ne 0$, then
\begin{align*}
\Omega_{\gamma} & = \sum_{s:\T^m_k(s) = \gamma}\sum_{x\in\F_{2^m}^*}(-1)^{\T_m(ux) \oplus \T_m(vsx^{2^k+1})} \\
& = \left\{\begin{array}{r@{\quad:\quad}l}
2^m-2^{m-k} & \gamma = \T_k^m(uv^d)^2,\\
-2^{m-k} & \mbox{otherwise.}
\end{array}\right.
\end{align*}
%where $\mu$ is the unique element in $\F_{2^m}$ for which $\mu^{2^k+1} = v^{-1}$.
\item[(ii)] If $u\ne 0$, then
\begin{align*}
\Upsilon_{\gamma} & = \sum_{s:\T^m_k(s) = \gamma}\sum_{x\in\F_{2^m}^*}(-1)^{\T_m(usx^{-d}) \oplus \T_m(vx)} \\
& = \left\{\begin{array}{r@{\quad:\quad}l}
2^m-2^{m-k} & \gamma^2 = \T_k^m(vu^e),\\
-2^{m-k} & \mbox{otherwise.}
\end{array}\right.
\end{align*}
%where $\mu$ is the unique element in $\F_{2^m}$ for which $\mu^{2^k+1} = v^{-1}$.
\end{itemize}
\end{lemma}
{\it Proof.}
(i) Let $s_\gamma\in\F_{2^m}$ such that $\T^m_k(s_\gamma) = \gamma$. Then $\{s\in\F_{2^m}\,:\,\T^m_k(s) = \gamma\} = \{s_\gamma \oplus s \,:\,\T^m_k(s) = 0\}$.
Denoting the subspace $\{s\in\F_{2^m} \,:\,\T^m_k(s) = 0\}$ by $\Lambda_0$ we obtain
\begin{align*}
\Omega_\gamma & = \sum_{s\in\Lambda_0}\sum_{x\in\F_{2^m}^*}(-1)^{\T_m(ux) \oplus \T_m(vs_\gamma x^{2^k+1} \oplus vsx^{2^k+1})} \\
& = \sum_{x\in\F_{2^m}^*}(-1)^{\T_m(ux \oplus vs_\gamma x^{2^k+1})}\sum_{s\in\Lambda_0}(-1)^{\T_m(vsx^{2^k+1})}.
\end{align*}
Observe that $\sum_{s\in\Lambda_0}(-1)^{\T_m(vsx^{2^k+1})} = 0$ if $vx^{2^k+1} \not\in\Lambda_0^\perp$, the orthogonal complement of $\Lambda_0$, and
$\sum_{s\in\Lambda_0}(-1)^{\T_m(vsx^{2^k+1})} = 2^{m-k}$ if $vx^{2^k+1} \in\Lambda_0^\perp$. Since for $z\in\F_{2^k}$ and $s\in\Lambda_0$ we have
$\T_m(zs) = \T_k(z\T_k^m(s)) = 0$, and since $\Lambda_0^{\perp}$ has dimension $k$, we have $\Lambda_0^{\perp} = \F_{2^k}$. Therefore,
\begin{equation}
\label{Og1}
\Omega_{\gamma} = 2^{m-k}\sum_{x\in\F_{2^m}^*:\atop vx^{2^k+1}\in\F_{2^k}^*}(-1)^{\T_m(ux\oplus s_\gamma vx^{2^k+1})}.
\end{equation}
As $\gcd(2^m-1,2^k+1) = 1$, there exists a unique element $\mu\in\F_{2^m}$ with $\mu^{2^k+1} = v^{-1}$, i.e., $\mu^e = v$, hence $\mu = v^d$.
Then $vx^{2^k+1} \in\F_{2^k}$ if and only if
$x^{2^k+1}\in\mu^{2^k+1}\F_{2^k} = \mu^{2^k+1}\F_{2^k}^{2^k+1}$. Hence $(\ref{Og1})$ reduces to
\[ \Omega_{\gamma} = 2^{m-k}\sum_{z\in\F_{2^k}^*}(-1)^{\T_m(u\mu z\oplus s_\gamma v\mu^{2^k+1}z^{2^k+1})} =
2^{m-k}\sum_{z\in\F_{2^k}^*}(-1)^{\T_m(u\mu z\oplus s_\gamma z^2)}. \]
Let $\nu = \T_k^m(u\mu) = \T_k^m(uv^d)$, %and $\tau\in\F_{2^k}$ be the (unique) square root of $\gamma$, then
then
\begin{align*}
\T_m(u\mu z\oplus s_\gamma z^2) & =  \T_k(z\T_k^m(u\mu)) \oplus \T_k(z^2\T_k^m(s_{\gamma})) = \T_k(\nu^2 z^2) \oplus \T_k(\gamma z^2) \\
& = \T_k((\nu^2 \oplus \gamma)z^2).
\end{align*}
Consequently,
\[ \Omega_{\gamma} = 2^{m-k}\sum_{z\in\F_{2^k}^*}(-1)^{\T_k((\nu^2 \oplus \gamma)z^2)} = \left\{\begin{array}{r@{\quad:\quad}l}
2^m-2^{m-k} & \nu^2\oplus\gamma = 0, \\
-2^{m-k} & \mbox{otherwise.}
\end{array}\right. \]
(ii) For $\Upsilon_{\gamma}$ with the same arguments we obtain
\[ \Upsilon_{\gamma} = 2^{m-k}\sum_{x\in\F_{2^m}^*:\atop ux^{-d}\in\F_{2^k}^*}(-1)^{\T_m(us_\gamma x^{-d}\oplus vx)}. \]
If $\mu\in\F_{2^m}$ is the unique element such that $\mu^d = u$, i.e., $\mu = u^e$, we get
\[ \Upsilon_{\gamma} = 2^{m-k}\sum_{z\in\F_{2^k}^*}(-1)^{\T_m(s_\gamma z^{-d} \oplus v\mu z)} = 2^{m-k}\sum_{z\in\F_{2^k}^*}(-1)^{\T_k(\gamma z^{-d}) \oplus \T_k(z\T_k^m(v\mu))}. \]
Changing the order of summation with $z\rightarrow z^{2^k+1} = z^{-e}$ we have
\begin{align*}
\Upsilon_{\gamma} & = 2^{m-k}\sum_{z\in\F_{2^k}^*}(-1)^{\T_k(\gamma (z^{-e})^{-d}) \oplus \T_k(z^{2^k+1}\T_k^m(v\mu))} \\
& = 2^{m-k}\sum_{z\in\F_{2^k}^*}(-1)^{\T_k(\gamma z) \oplus \T_k(z^2\T_k^m(v\mu))} = 2^{m-k}\sum_{z\in\F_{2^k}^*}(-1)^{\T_k((\gamma^2 \oplus \T_k^m(v\mu))z^2)},
\end{align*}
which yields the claimed result for $\Upsilon_\gamma$.
\hfill$\Box$
\begin{theorem}
\label{partithe}
Let $m,k$ be integers such that $k$ divides $m$ and $\gcd(2^m-1,2^k+1) = 1$, and let $\pi(i) = \gamma_i$ be a one-to-one map from $\Z_{2^k}$ to $\F_{2^k}$.
%w.l.o.g. let $\gamma_0 = 0$ (hence w.l.o.g., $\pi(0) =0$).
Define functions $f_A, f_B:\F_{2^m}\times\F_{2^m}\rightarrow\Z_{2^k}$ as follows:
\begin{itemize}
\item[-] If $(x,y)\in\mathcal{A}(\gamma_i)$ then $f_A(x,y) = i$, and, w.l.o.g., $f_A(0,y) = 0$ for all $y\in\F_{2^m}$;
\item[-] If $(x,y)\in\mathcal{B}(\gamma_i)$ then $f_B(x,y) = i$, and, w.l.o.g., $f(x,0) = 0$ for all $x\in\F_{2^m}$.
\end{itemize}
Then $f_A, f_B$ are $\Z_{2^k}$-bent functions.
\end{theorem}
{\it Proof.} For $u,v\in \F_{2^m}$ and $0\le t\le k-1$, for the function $f_A$ we have
\begin{align*}
\mathcal{H}_{f_A}(2^t,(u,v)) & = \sum_{x,y\in\F_{2^m}}\zeta_{2^k}^{2^tf_A(x,y)}(-1)^{\T_m(uv\oplus vy)} \\
& = \sum_{i=0}^{2^k-1}\sum_{(x,y)\in\mathcal{A}(\gamma_i)}\zeta_{2^k}^{2^ti}(-1)^{\T_m(ux \oplus vy)} + \sum_{y\in\F_2^m}(-1)^{\T_m(vy)} \\
& = \sum_{i=0}^{2^k-1}\zeta_{2^k}^{2^ti}\sum_{\T_k^m(s)=\gamma_i}\sum_{x\in\F_{2^m}^*}(-1)^{\T_m(ux\oplus vsx^{2^k+1})} + \sum_{y\in\F_{2^m}}(-1)^{\T_m(vy)}.
\end{align*}
If $v=0$, then
\[ \mathcal{H}_{f_A}(2^t,(u,v)) = 2^{m-k}\sum_{x\in\F_{2^m}^*}(-1)^{\T_m(ux)}\sum_{i=0}^{2^k-1}\zeta_{2^k}^{2^ti} + \sum_{y\in\F_{2^m}}1 = 2^m. \]
Let $v\ne 0$. Then with $\sum_{(x,y)\in\mathcal{A}(\gamma_j)}(-1)^{\T_m(ux\oplus vy)} = \Omega_\gamma$, and applying Lemma \ref{aboutO}, we obtain
\begin{align}
\label{iota1}
\nonumber
\mathcal{H}_{f_A}(2^t,(u,v)) & = \sum_{i=0}^{2^k-1}\zeta_{2^k}^{2^ti}\Omega_{\gamma_i} + \sum_{y\in\F_{2^m}}(-1)^{\T_m(vy)} =
-2^{m-k}\sum_{ci=0}^{2^k-1}\zeta_{2^k}^{2^ti} + 2^m\zeta_{2^k}^{2^t\iota} \\
& = 2^m\zeta_{2^k}^{2^t\iota},
\end{align}
if $\gamma_\iota = \T_k^m(u\mu)^2$, where $\mu^{2^k+1} = v^{-1}$. \\[.3em]
In the same way, for $f_B$ we obtain $\mathcal{H}_{f_B}(2^t,(0,v)) = 2^m$ (as w.l.o.g. $f_B(x,0) = 0$), and
\begin{equation}
\label{iota2}
\mathcal{H}_{f_B}(2^t,(u,v)) = 2^m\zeta_{2^k}^{2^t\iota}
\end{equation}
if $\T_k^m(u^ev) =\gamma^2_\iota$.
\hfill$\Box$ \\[.5em]
From Theorem \ref{partithe} we follow the next theorem on the related Boolean bent functions, which, just as the spread bent functions, are exactly those
bent functions which are constant on the elements of a certain partition of $\V_n$.
\begin{theorem}
\label{PSap}
Let $m,k$ be integers such that $k$ divides $m$ and $\gcd(2^m-1,2^k+1) = 1$, let $e=2^m-2^k-2$ and $d$ such that $de\equiv 1\bmod 2^m-1$.
For $\gamma\in\F_{2^k}$ let
\[ \mathcal{A}(\gamma) = \bigcup_{s\in\F_{2^m}:\T_k^m(s) = \gamma}\{(x,sx^{2^k+1})\,:x\in\F_{2^m}^*\}, U = \{(0,y)\,:\,y\in\F_{2^m}\}, \]
and
\[ \mathcal{B}(\gamma) = \bigcup_{s\in\F_{2^m}:\T_k^m(s) = \gamma}\{(x^{-d}s,x)\,:x\in\F_{2^m}^*\}, V = \{(x,0)\,:\,x\in\F_{2^m}\}. \]
\begin{itemize}
\item[I.] Every Boolean function whose support is the union of $2^{k-1}$ of the sets $\mathcal{A}(\gamma)$ is a bent function.
Likewise, their complements, i.e., the Boolean functions with $U$ and $2^{k-1}$ of the sets $\mathcal{A}(\gamma)$ as their support, are bent.
\item[II.] Every Boolean function whose support is the union of $2^{k-1}$ of the sets $\mathcal{B}(\gamma)$ is a bent function.
Likewise the Boolean functions with $V$ and $2^{k-1}$ of the sets $\mathcal{B}(\gamma)$ as their support, are bent.
\end{itemize}
The duals of the bent functions of the class in I are in the class in II (and vice versa).
\end{theorem}
{\it Proof.}
The bentness of the Boolean functions in I and II follows immediately with $(\ref{iota1})$ and $(\ref{iota2})$ for $t = k-1$. Note that the Boolean function
$2^{k-1}f$ has as support all $\mathcal{A}(\gamma)$ respectively $\mathcal{B}(\gamma)$ which $f$ maps to an odd $i$. Moreover, we explicitly see the
dual function of a bent function of I from $(\ref{iota1})$ (for $v\ne 0$)
\[ \mathcal{H}_f(2^{k-1},(u,v)) = 2^m(-1)^\iota \]
if $\gamma_\iota = \T_k^m(uv^d)^2$. Again observing that $\T_k^m(uv^d)^2$ is constant $\gamma$ for all $(u,v) \in \mathcal{B}(\gamma)$, we infer that
the dual is in class II. Note that clearly, as the dual is bent, we must have exactly $2^{k-1}$ of the sets $\mathcal{B}(\gamma)$ in the support
(additionally $V$ may be in the support). % - in Theorem w.l.o.g $v=0$ then $2^m$).
\hfill$\Box$
\begin{remark}
The Maiorana-McFarland functions $\T_m(\alpha^{-e} x^ey)$, $\alpha\in\F_{2^k}^*$, and the functions $g(x)$ in Corollary \ref{Carletsec} certainly belong to the
functions of class I, their duals $\T_m(\alpha xy^d)$ and $g^*$ in Corollary \ref{Carletsec} belong to class II. As easily seen, in general these bent functions and their duals have different algebraic degree, hence in general they are not 
EA-equivalent. These spread-like partitions, which also give rise to $\Z_{2^k}$-bent functions and are the main
object of our interest, yield a large quantity of bent functions of class I respectively class II. If some given function
in class I or class II can be obtained from any of the numerous known primary or secondary bent function constructions
is as usual very difficult to answer. An arbitrary function in class I or class II does not necessarily belong to the  
completed Maiorana-McFarland class as Example 6.3.16 in \cite{dillon} shows (special case $k=m$, where our partition
reduces to the Desarguesian spread - see the remark below).
\end{remark}

\begin{remark}
\begin{itemize}
\item[(i)] In the special case $k=m$, the set $\mathcal{S} = \{U,\mathcal{A}(\gamma)\;:\;\gamma\in\F_{2^m}\}$ reduces to the Desarguesian spread, and $f$ in
Theorem \ref{partithe} is a spread function on the complete Desarguesian spread, described as in Section \ref{prel}. Corollary \ref{PSap} describes then the PS$_{ap}^-$ and
PS$_{ap}^+$ functions, cf. \cite{dillon}. Hence we may see the $\Z_{2^k}$-bent functions in Theorem \ref{partithe}, and the Boolean bent functions in Corollary \ref{PSap} as generalizations
of the Desarguesian spread bent functions.
\item[(ii)] As for the spread functions in Section \ref{prel}, also the proof of Theorem \ref{partithe}, holds not only for functions from $\F_{2^m}\times\F_{2^m}$ to
$\Z_{2^k}$, but for functions from $\F_{2^m}\times\F_{2^m}$ into any abelian group of order $2^k$. The bentness is a property of the partition of $\F_{2^m}\times\F_{2^m}$.
In particular, also many more vectorial bent functions in dimension $k$ are obtained. As easily observed, these vectorial bent functions from $\F_{2^m}\times\F_{2^m}$ to
$\F_{2^k}$ belong to the class of vectorial dual-bent functions, where the duals for the functions that are constant on $\Gamma_1$ are functions that are constant on $\Gamma_2$.
%\item[(iii)] Bent functions respectively vectorial bent functions we get are not?? MMF, that would be good.
\end{itemize}
\end{remark}
For a spread $\mathcal{S} = \{U_j\;:\;j=0,1,\ldots,2^m+1\}$ of $\V_n$ denote (for a fixed inner product) by $\mathcal{S}^\perp$ the spread
$\mathcal{S}^\perp = \{U_j^\perp\;:\;j=0,1,\ldots,2^m+1\}$.

As is well known, the duals of the Boolean partial spread bent functions from the spread $\mathcal{S}$ are the Boolean spread bent functions from $\mathcal{S}^\perp$.
Moreover, a vectorial spread bent function from $\mathcal{S}$ is vectorial dual-bent, with a dual obtained from $\mathcal{S}^\perp$  (see \cite[Theorem 3]{cmp18}).
Hence also with this respect the functions obtained from the partitions $\Gamma_1$ and $\Gamma_2$ behave like spread bent functions. In fact, for the special case
that $e\equiv 2^v\bmod 2^m-1$, in which case $\Gamma_1, \Gamma_2$ are spreads, we see that $\Gamma_2^\perp = \Gamma_1$.
For the Desarguesian spread (more general for every symplectic spread), the Boolean spread bent functions and their duals are EA-equivalent, see the discussion on weak
self-duality in \cite{cmp18}. As observed above, with this respect the (vectorial) bent functions obtained from $\Gamma_1,\Gamma_2$ in general behave different.

We close this section with some more remarks on equivalence of partitions. We here call two partitions of $\V_n$
equivalent, if one is obtained from the other with a linear coordinate transformation on $\V_n$. 
Clearly, as for the spreads, the partitions $\Gamma_1$ and $\Gamma_2$ represent then a whole equivalence class
of partitions. Let $\mathcal{L}$ be a linear permutation of $\F_{2^m}\times\F_{2^m}$ and let $\mathcal{L}(\Gamma_i)$, $i=1,2$, be the partition of $\F_{2^m}\times\F_{2^m}$
obtained via the coordinate transformation described by $\mathcal{L}$. Then every Boolean function $g$ which is constant on the elements of $\mathcal{L}(\Gamma_i)$  and
has a support of cardinality $2^{2m-1}-2^{m-1}$ or $2^{2m-1}+2^{m-1}$ is of the form $g(x) = f(\mathcal{L}^{-1}(x))$ for a bent function in the class I respectively II.
Hence $g$ is a bent function affine equivalent to $f$. Obviously two such partitions which are equivalent via a coordinate transformation yield the same sets of bent functions (in the sense of EA-equivalence for Boolean functions).

Though we observed that we have functions $f_1$ from class I and functions $f_2$ from class II which are not 
EA-equivalent, the partitions $\Gamma_1$ and $\Gamma_2$ may still
be equivalent (note that in general also the set of bent functions for a given spread contains several EA-equivalence classes). We remark that for the Desarguesian spread (and in general
for symplectic spreads) $\mathcal{S}$ and $\mathcal{S}^\perp$ are equivalent. Numerically we confirmed that in general $\Gamma_1$ and $\Gamma_2$ are not equivalent.
Hence with respect to equivalence, the partitions $\Gamma_1$, $\Gamma_2$ behave different than $\mathcal{S}$ and $\mathcal{S}^\perp$ for a Desarguesian spread, the special
case when $e\equiv 2^v\bmod 2^m-1$. \\[.5em]

\noindent
{\bf Acknowledgement.} \\[.5em]
The authors are supported by the FWF Project P 30966.


\begin{thebibliography}{99}

\bibitem{c06} C. Carlet, On bent and highly non-linear balanced/resilient functions and their algebraic immunities,
in: M.P.C. Fossorier et al. (Eds.), AAECC, Lecture Notes in Computer Science 3857, pp. 1--28, Springer-Verlag, New York, 2006.

\bibitem{c11} C. Carlet, Relating three nonlinearity parameters of vectorial functions and building APN functions from bent functions.
Des. Codes Cryptogr. 59 (2011), 89--109.

\bibitem{cmp18} A. \c Ce\c smelio\u glu, W. Meidl, A. Pott, Vectorial bent functions and their duals. Linear Algebra Appl. 548 (2018), 305--320.

\bibitem{dillon} J.F. Dillon, Elementary Hadamard difference sets, Ph.D. dissertation, University of Maryland, 1974.

\bibitem{hmp} S. Hod\v zi\'c, W. Meidl, E. Pasalic, Full characterization of generalized bent functions as (semi)-bent spaces,
their dual, and the Gray image. IEEE Trans. Inform. Theory 64 (2018), 5432--5440.

\bibitem{ll} P. Lisonek, H.Y. Lu, Bent functions on partial spreads. Des. Codes Cryptogr. 73 (2014), 209--216.

\bibitem{mms} T. Martinsen, W. Meidl, P. Stanica, Generalized bent functions and their Gray images. In: Arithmetic of finite fields,
Lecture Notes in Comput. Sci., 10064, pp. 160--173, Springer, Cham, 2016.

\bibitem{mms1} T. Martinsen, W. Meidl, P. Stanica, Partial spread and vectorial generalized bent functions. Des. Codes Cryptogr. 85 (2017), 1--13.

\bibitem{m18} W. Meidl, A secondary construction of bent functions, octal gbent functions and their duals. Math. Comput. Simulation 143 (2018), 57--64.

\bibitem{mp} W. Meidl, A. Pott, Generalized bent functions into $\Z_{p^k}$ from the partial spread and the Maiorana-McFarland class,  Cryptogr. Commun. 11 (2019), 1233--1245.

\bibitem{mtqwwf} S. Mesnager, C. Tang, Y. Qi, L. Wang, B. Wu, K. Feng, Further results on generalized bent functions and their
complete characterization. IEEE Trans. Inform. Theory 64 (2018), 5441--5452.

\bibitem{mes} S. Mesnager, Several new infinite families of bent functions and their duals. IEEE Trans. Inform. Theory 60 (2014), no. 7, 4397--4407.

\bibitem{mesbook} S. Mesnager, Bent functions. Fundamentals and results. Springer, 2016.

\bibitem{n} K. Nyberg, Perfect nonlinear S-boxes, In: Advances in cryptology--EUROCRYPT '91 (Brighton, 1991),
Lecture Notes in Comput. Sci., 547, pp. 378--386, Springer, Berlin, 1991.

\bibitem{p04} A. Pott, Nonlinear functions in abelian groups and relative difference sets. Discrete Applied Mathematics 138 (2004), 177--193.

\bibitem{p96} A. Pott, A survey on relative difference sets. Groups, difference sets, and the Monster, In: Ohio State Univ. Math. Res. Inst. Publ., 4,
pp. 195--232, de Gruyter, Berlin, 1996.

\bibitem{sb} B. Schmidt, On $(p^a,p^b,p^a,p^{a-b})$-relative difference sets. J. Algebraic Combin. 6 (1997), 279--297.

\bibitem{txqf} C. Tang, C. Xiang, Y. Qi, K. Feng, Complete characterization of generalized bent and $2^k$-bent Boolean functions.
IEEE Trans. Inform. Theory 63 (2017), 4668--4674.

\end{thebibliography}
\end{document}